\documentclass{amsart}
\usepackage{multirow}
\usepackage{amssymb,  latexsym, ulem}

\usepackage{cancel}

\usepackage[all]{xy}

\newtheorem{theorem}{Theorem}[section]
\newtheorem{lemma}[theorem]{Lemma}
\newtheorem{sublemma}[theorem]{Sublemma}
\newtheorem{corollary}[theorem]{Corollary}
\newtheorem{proposition}[theorem]{Proposition}

\theoremstyle{definition}
\newtheorem{definition}[theorem]{Definition}

\theoremstyle{remark}
\newtheorem{remark}[theorem]{Remark}

\newcommand{\lem}[2]{\begin{lemma}\label{#1} #2\end{lemma}}
\newcommand{\lemc}[3]{\begin{lemma}[\mbox{#2}]\label{#1} #3\end{lemma}}

\newcommand{\prop}[2]{\begin{proposition}\label{#1}#2\end{proposition}}

\newcommand{\thm}[2]{\begin{theorem}\label{#1}#2\end{theorem}}

\newcommand{\cor}[2]{\begin{corollary}\label{#1}#2\end{corollary}}

\newcommand{\al}{\alpha}
\newcommand{\be}{\beta}
\newcommand{\ga}{\gamma}
\newcommand{\Ga}{{BP_*(BP)}}
\newcommand{\de}{\delta}
\newcommand{\De}{{\Delta}}
\newcommand{\ep}{\varepsilon}
\newcommand{\ze}{\zeta}
\newcommand{\et}{\eta}

\newcommand{\io}{\iota}
\newcommand{\ka}{\kappa}

\newcommand{\rh}{\rho}
\newcommand{\si}{\sigma}
\newcommand{\Si}{{\Sigma}}

\newcommand{\om}{\omega}

\newcommand{\Om}{{\Omega}}
\newcommand{\barr}{\begin{array}}
\newcommand{\ear}{\end{array}}
\newcommand{\bsk}{\begin{array}{rcl}}
\newcommand{\esk}{\end{array}}
\newcommand{\skh}{\begin{eqnarray*}}
\newcommand{\sko}{\end{eqnarray*}}
\newcommand{\bcss}{\begin{cases}}
\newcommand{\ecss}{\end{cases}}

\newcommand{\sk}[1]{$$\begin{array}{rl}#1\end{array}$$}
\newcommand{\skr}[2]{$$\begin{array}{rl}#1\end{array}\lnr{#2}$$}
\newcommand{\skc}[3]{$$\begin{array}{rl}#1\end{array}\leqno{(\thetheorem.#2)\label{#3}}$$}
\newcommand{\cgs}{\cg & \hsp{-.1in} }
\newcommand{\urcgs}[1]{\underset{\ref{#1}}{\cg} & \hsp{-.1in} }
\newcommand{\ukcgs}[1]{\underset{\kko{#1}}{\cg} & \hsp{-.1in} }
\newcommand{\usscgs}[1]{\underset{\scriptsize \begin{matrix}#1\end{matrix}}{\cg} & \hsp{-.1in} }
\newcommand{\eqs}{= & \hsp{-.1in} }

\newcommand{\inds}[1]{ & \hsp{-.#1in} }

\newcommand{\AL}[1]{\begin{align*}#1\end{align*}}
\newcommand{\ALr}[2]{$$\begin{aligned}#1\end{aligned}\lnr{#2}$$}

\newcommand{\AR}[2]{$$\begin{array}{#1}#2\end{array}$$}
\newcommand{\ARr}[3]{$$\begin{array}{#1}#2\end{array}\lnr{#3}$$}
\newcommand{\aR}[2]{$\begin{array}{#1}#2\end{array}$}
\newcommand{\naka}[1]{\begin{center}#1\end{center}}

\newcommand{\LE}[1]{\lefteqn{#1}}
\newcommand{\cass}[1]{\begin{cases}#1\ecss}
\newcommand{\ak}{\quad}

\renewcommand{\Xy}[2]{$$\xymatrix@=#1pt{#2}$$}
\newcommand{\Xyr}[3]{$$\vcenter{\xymatrix@=#1pt{#2}}\lnr{#3}$$}
\newcommand{\Xyrs}[4]{$$\vcenter{\xymatrix@=#1pt{#2}}\lnrs{#4}{#3}$$}
\newcommand{\xR}[2]{$\xymatrix@R#1pt{#2}$}
\newcommand{\xyR}[2]{$$\xymatrix@R#1pt{#2}$$}
\newcommand{\xyRr}[3]{$$\vcenter{\xymatrix@R#1pt{#2}}\lnr{#3}$$}
\newcommand{\xyRrs}[4]{$$\vcenter{\xymatrix@R#1pt{#2}}\lnrs{#4}{#3}$$}

\newcommand{\wt}[1]{\widetilde{#1}}

\newcommand{\N}{\mathbb N}
\newcommand{\Nt}[1]{{\mathbb N}^{(#1)}}
\newcommand{\Z}{\mathbb Z}

\newcommand{\Zpp}{{\mathbb Z}^{(p)}}

\def\o+{\oplus}
\newcommand{\Op}{\bigoplus}
\newcommand{\x}{\times}
\newcommand{\ox}{\otimes}

\newcommand{\e}{{\rm Ext}}
\newcommand{\lnb}{\refstepcounter{theorem}\leqno(\thetheorem)}
\newcommand{\lnr}[1]{\lnb\label{#1}}
\newcommand{\nr}{\refstepcounter{theorem}\thetheorem}
\newcommand{\kko}[1]{(\ref{#1})}
\newcommand{\mbx}[1]{\quad\mbox{#1}\quad}
\newcommand{\cf}{{\it cf.}\ }

\newcommand{\Ker}{{\rm Ker}\ }

\newcommand{\im}{{\rm Im}\, }

\newcommand{\ds}{\displaystyle }
\newcommand{\cg}{\equiv}

\newcommand{\qand}{\mbx{and}}
\newcommand{\eR}{\et_R}
\newcommand{\Lt}{\left}
\newcommand{\Rt}{\right}
\newcommand{\LR}[1]{\Lt(#1\Rt)}
\newcommand{\cln}{\colon}
\renewcommand{\O}[1]{\overline{#1}}
\newcommand{\C}[2]{{#1\choose #2}}

\newcommand{\mx}[1]{\begin{matrix}#1\end{matrix}}

\newcommand{\da}{\mathfrak{a}}
\newcommand{\db}{\mathfrak{b}}
\newcommand{\dc}{\mathfrak{c}}
\newcommand{\dg}{\mathfrak{g}}

\newcommand{\dk}{\mathfrak{k}}
\newcommand{\dx}{\mathfrak{x}}
\newcommand{\Dp}{\mathfrak{p}}

\newcommand{\dz}{\mathfrak{z}}

\newcommand{\xar}{\xrightarrow}
\newcommand{\U}[2]{\uline{#1}_{#2}}

\newcommand{\cS}{{\mathcal S_{(p)}}}

\newcommand{\wh}[1]{\widehat{#1}}

\newcommand{\p}{\noindent {\it Proof.} }
\newcommand{\q}{\hfill $\qed$ \medskip}

\newcommand{\A}{{BP_*}}

\renewcommand{\b}[2]{b_{#1,#2}}

\newcommand{\tp}[1]{t_1^{p^{#1}}}
\newcommand{\pt}[1]{^{p^{#1}}}

\newcommand{\usk}[1]{\underset{\kko{#1}}}
\newcommand{\usr}[1]{\underset{\ref{#1}}}
\newcommand{\uss}[1]{\underset{\scriptsize \begin{matrix}#1\end{matrix}}}
\newcommand{\hsp}{\hspace}
\newcommand{\vsp}{\vspace}

\newcommand{\vm}{v_{n-1}}

\newcommand{\hc}[1]{\left[#1\right]}
\newcommand{\B}[1]{E(#1)_*}
\newcommand{\Sg}[1]{E(#1)_*(E(#1))}
\newcommand{\kt}[2]{(t_1^{#1})_{#2}}
\newcommand{\kh}[2]{(h_{#1})_{#2}}
\newcommand{\UW}{\uwave}

\newcommand{\Blk}{\Big(}
\newcommand{\Brk}{\Big)}

\begin{document}

\title[On products of beta and gamma elements]{On products of beta and gamma elements in the homotpy of the first Smith-Toda spectrum }
\author{ Katsumi Shimomura}
\address{Department of Mathematics, Faculty of Science, Kochi University, Kochi,
780-8520, Japan}
\email{katsumi@kochi-u.ac.jp}
\author{Mao-no-suke Shimomura}
\address{Department of Mathematics, Faculty of Science and Technology, Kochi University, Kochi,
780-8520, Japan}
\email{b21m6g02@s.kochi-u.ac.jp}


\begin{abstract}
In this paper, we determine the first cohomology  of the monochromatic comodule $M^1_2$ at an odd prime, and apply the results to show non-trivialities of some products of beta and gamma elements in the homotopy groups of the Smith-Toda spectrum $V(1)$. The cohomology for  a prime greater than three was determined by the first author \cite{sV1}. Here, we verify them and determine the cohomology at the prime 3 by elementary calculation. The cohomology will be a stepping stone for computing the cohomology of the monochromatic comodule $M^3_0$, which we hope to determine for a long time.
\end{abstract}
\maketitle

\section{Introduction}

Let $p$ be an odd prime number, and $\cS$ denote the stable homotopy category of $p$-local spectra.
Let $S\in \cS$ denote the sphere spectrum.
Then, the mod $p$ Moore spectrum $M$ and the first Smith-Toda spectrum $V(1)$ are given by the cofiber sequences
$$
S\xar{p} S\xar iM\xar j\Si S\qand \Si^qM\xar{\al} M\xar{i_1}V(1)\xar{j_1}\Si^{q+1}M.\lnr{cof1}
$$
Here, $p\in \pi_0(S)\cong \Z_{(p)}$, and $\al\in [M,M]_q$ denotes the Adams map.
Hereafter, we put
$$
q=2p-2\in \Z.
$$
In order to study the homotopy groups $\pi_*(X)$ of a spectrum $X$, we adopt the Adams-Novikov spectral sequence
$$
E_2^{s,t}(X)=H^{s,t}BP_*(X)\Longrightarrow \pi_{t-s}(X). \lnr{ANSS}
$$
Hereafter, we abbreviate as 
$$
H^{s,t}M=\e_\Ga^{s,t}(\A,M)
$$
for a $\Ga$-comodule $M$ over the Hopf algebroid
$$
(BP_*, BP_*(BP))=(\Z_{(p)}[v_1,v_2,\dots], BP_*[t_1,t_2,\dots]) \lnr{Hal}
$$
based on the Brown-Peterson spectrum $BP\in \cS$.
We note that $v_i$'s are Hazewinkel's generators and the degrees of $v_i$ and $t_i$ are $|v_i|=2p^i-2=|t_i|$.

Let 
$$
I_n=(p,v_1,\dots,v_{n-1}) \qand J_j=(p,v_1,v_2^j) \lnr{ideals}
$$ 
($v_0=p$) denote the invariant ideals of $BP_*$.
Since $BP_*(\al)=v_1$, the cofiber sequences \kko{cof1} induce the short exact sequences
\ARr c{
0\to \A\xar{p}\A\xar{i_*} \A/I_1\to 0\qand\\
 0\to \A/I_1\xar{v_1}\A/I_1\xar{(i_1)_*} \A/I_2\to 0
}{ses0}
along with the isomorphisms
$$
\A(S)=\A,\ak \A(M)=\A/I_1,\qand \A(V(1))=\A/I_2.
$$
Furthermore, we have a short exact sequence
$$
0\to \A/I_2\xar{v_2^j}\A/I_2\xar{\O i_j}\A/J_j\to 0 \lnr{ses2}
$$
for $j\ge 1$.
We denote by $\de_0\cln H^s\A/I_1\to H^{s+1}\A$, $\de_1\cln H^s\A/I_2\to H^{s+1}\A/I_1$ and $\O\de_j\cln H^s\A/J_j\to H^{s+1}\A/I_2$, the connecting homomorphisms associated to the short exact sequences \kko{ses0} and \kko{ses2}.  
We define the Greek letter elements by:
\AR{rlll}{
\O\be_s'&\hsp{-.1in}=\de_1(v_2^s)&\in E_2^1(M)=H^1\A/I_1&\mbox{for $v_2^s\in H^0\A/I_2$,}\\
\O\be_s&\hsp{-.1in}=\de_0\de_1(v_2^s)&\in E_2^2(S)=H^2\A&\mbox{for $v_2^s\in H^0\A/I_2$, and}\\
\O\ga_{s/j}''&\hsp{-.1in}=\O\de_j(v_3^s)&\in E_2^1(V(1))=H^1\A/I_2&\mbox{for $v_3^s\in H^0\A/J_j$,}
}
and $\O\ga''_s=\O\ga_{s/1}''\in E_2^1(V(1))$. We notice that $1\le j\le p^n$ if ${p^n| s}$.

Let $\Z$ and $\N$ denote the set of all integers and  its subset consisting of all non-negative integers, respectively.
We denote by $\Z^{(p)}(=\Z\setminus p\Z)$ and $\N^{(p)}(=\N\setminus p\N)$  the set of the integers prime to $p$, and decompose $\Zpp$ into the three summands:
\skr{
\Zpp\eqs \Z_0\coprod \Z_1 \coprod \Z_2, \mbx{for}\\
\Z_0\eqs\{s\in\Zpp\mid  p\nmid (s+1)\},\ak
\Z_1= \{s\in\Zpp\mid  p^2| (s+1)\},\qand\\
\Z_2\eqs \{s\in\Zpp\mid  p| (s+1)\ \mbox{and}\ p^2\nmid (s+1)\}.
}{Zi}
We consider subsets of $\N$: 
\sk{
2\N_{>0}\eqs \{s\in \N\mid s \mbox{ is even}\ge 2\}, \qquad
\O{2\N}= \{s\in \N\mid s \mbox{ is odd}\},\\
\N_{1}\eqs \{s\in\N^{(p)}\mid  p^2\nmid (s+p+1),\mbox{ or } p^3\mid (s+p+1)\}, \qand\\
\N_2\eqs \{s\in\Nt p\mid p\nmid (s+2),\mbox{ or } p^3\mid (s+2)(s+2+p)\}.
}
Furthermore,  we put $\Z_i^+=\Z_i\cap \N$ for $i=0,1,2$.
We introduce the subsets $U$, $U_1$ and $U_2$ of $\Nt p\x \N$ given by
\AR l{
U_1= (\Nt p\x 2\N)\cup (\Z_0^+\x\N), \\
U'_1= (\Nt 3\x \{0\})\cup (\N_1\x 2\N_{>0})\cup ((\Z_0^+\cap \N_2)\x \N)\cup (\Z_0^+\x \{1\})\\
U_2= (\N_1\x 2\N)\cup ((\Z_0^+\cap\N_2)\cup \Z_1^+)\x \N)\cup (\Nt p\x \{1\})\qand \\
U_2'= (\N_1\x\{0\})\cup (\Nt 3\x (\{1\}\cup 2\N_{>0})\cup ((\Z_0^+\cup \Z_1^+)\x \N)
}

Our main result is the following:

\thm{main1}{Let $p$ be an odd prime.
In the Adams-Novikov $E_2$-term for computing $\pi_*(V(1))$, $\O\be_1$ and $\O\be_2$ act on the gamma elements $\O\ga''_{sp^r/j}$ 
$((s,r)\in\Nt p\x \N$ and $1\le j\le p^r)$ by: 
\sk{
\O\ga''_{sp^r/j}\O\be_1\ne 0 &\mbox{for  $(s,r)\in U_1$ if $p\ge 5$, and
for $(s,r)\in U_1'$ if $p=3$, }\\
\O\ga''_{sp^r/j}\O\be_2\ne 0&\mbox{for $(s,r)\in U_2$ if $p\ge 5$, 
 and for $(s,r)\in U_2'$ if $p=3$.} 
}
in $E_2^3(V(1))$.
}
We notice that there is a way to define $\ga''_{sp^r/j}$ for $j\le a_{r}$ ($a_r$ in \kko{ea}) so that $v_2^{j-1}\ga''_{sp^r/j}=\ga''_{sp^r}$, and the theorem holds for such extended gamma elements. 
We also notice that $\O \be_s\cg \C s2v_2^{s-2}\O\be_2+s(2-s)v_2^{s-1}\O \be_1$ mod $I_2$ (\cf \cite[Lemma 4.4]{os}), and so
$$
\O\ga''_{sp^r/j}\O\be_t=\C t2\O\ga''_{sp^r/j-t+2}\O\be_2+t(2-t)\O\ga''_{sp^r/j-t+1}\O\be_1.
$$
Thus, Theorem \ref{main1} implies  non-triviality of the products of $\O\ga''_{sp^r/j}$ and $\O\be_t$.

The Adams-Novikov differential $d_r=0$ if $q\nmid (r-1)$ by the sparseness of the spectral sequence \kko{ANSS}. This shows that the products in the theorem are not in the image of any differentials $d_r$.
It is well known that the elements $\O \be_1$ and $\O\be_2$ converge to the homotopy elements $\be_1$ and $\be_2\in \pi_*(S)$, respectively, in the spectral sequence \kko{ANSS} for $X=S$.

\cor{main2}{Let $p$ be an odd prime. If $\O\ga''_{sp^r/j}\in E_2^1(V(1))$ is a permanent cycle detecting $\ga''_{sp^r/j}\in \pi_*(V(1))$,
then, $\ga''_{sp^r/j}\be_i\ne 0$ $(i=1,2)$ in the homotopy groups $\pi_*(V(1))$ for $(s,r)$ given in Theorem \ref{main1}.
}

Toda \cite[Th.~1]{toda} and Oka \cite[Th.~4.2]{ok} showed that $\ga''_s$ and $\ga''_{sp/2}$ are permanent cycles for $p\ge 7$. 

\cor{main3}{Let $p\ge 7$ and $r$ and $s$ be integers with $(s,r)\in\Nt p\x \N$. Then,  in $\pi_*(V(1))$,
\AR l{
\ga''_{sp^r/j}\be_1\ne 0\mbx{if $r$ is even or $p\nmid (s+1)$, }\\
\ga''_{sp^{2r}/j}\be_2\ne 0\mbx{if $p^2\nmid (s+p+1)$ or $p^3|(s+p+1)$, }\\
\ga''_{sp^{2r+1}/j}\be_2\ne 0\mbx{for $r\ge 1$ if $p\nmid (s+1)(s+2)$, $p^2|(s+1)$ or $p^3|(s+2)(s+2+p)$.}
}
and $\ga''_{sp/j}\be_2\ne 0$, where
 $j=1,2$.
}
Theorem \ref{main1} follows from Theorem \ref{H1}, which states the structure of the first cohomology of the monochromatic comodule $M^1_2$.
The cohomology $H^1M^1_2$ was determined by the first author \cite{sV1} based on the computation in \cite{seR} at a prime $\ge 5$.  
In this paper, we determine the cohomology based on elementary calculation at an odd prime. 
The generators are explicitely given so that we will use the result easily in further computation.
This result will be a stepping stone for determining the long desired cohomology $H^*M^3_0$.

This paper is organized as follows:
In the next section, we state the main result, Theorem \ref{H1}, which gives the structure of  $H^1M^1_2$.
In section three, we prove Theorems \ref{H1} and \ref{main1} assuming Lemma \ref{key}, whose proof  will be given in the last section.
Section four is devoted to introducing some formulas, cochains and relations for the following sections.
We refine the elements $x_{3,i}$ given in \cite[(5.11)]{mrw} to define $x_i$, which induce the cochains $y_{s,i}$ and $y_{s,i}'\in\Om^1E(3)_*$
in section five.

\section{The structure of  $H^1M^1_2$}

In this section, we state the structure of $H^1M^1_2$ for an odd prime $p$ obtained in this paper.
The structure was given in \cite{sV1}, which was done for the prime $p\ge 5$.

We begin with defining the monochromatic $\Ga$-comodules $N^s_n$ and $M^s_n$ inductively by
$$
N^0_n=\A/I_n,\ak M^s_n=v_{s+n}^{-1}N^s_n
$$
for the ideal $I_n$ in \kko{ideals}
and the short exact sequence
$$
0\to N^s_n\xar {\io^s_n}M^s_n\xar{\ka^s_n}N^{s+1}_n\to 0
\lnr{ses:ch}
$$
(\cite[\S 3.~A.]{mrw}). Since $\A$ is a $\Ga$-comodule with structure map $\eR$, the right unit map of the Hopf algebroid $\Ga$,
these monochromatic comodules have the structure maps induced from  $\eR$.

Let $E(3)$ denote the third Johnson-Wilson spectrum, which yields a Hopf algebroid
$$
(E(3)_*, E(3)_*(E(3)))=(\Z_{(p)}[v_1,v_2, v_3,v_3^{-1}], E(3)_*\ox_{\A}\Ga\ox_{\A}E(3)_*).
$$
Its structure maps are
induced from the Hopf algebroid $(\A,\Ga)$ in \kko{Hal}. 
Since we have the Miller-Ravenel change of rings theorem 
$$
H^*M=\e_{\Ga}^*(\A,M)\cong \e_{\Sg3}^*(\B3,\B3\ox_\A M)
$$ for a $v_3$-local $\Ga$-comodule $M$ (\cite[Th.~3.10]{m-r}),
we denote
 the cohomology of an $\Sg3$-comodule $M$ also by 
$$
H^sM=\e_{\Sg3}^s(\B3,M).
$$
By virtue of the change of rings theorem, we denote simply by $M^s_n$ the $\Sg3$-comodule $E(3)_*\ox_{\A}M^s_n$. 
In this paper, we consider the Ext group as the cohomology group of the cobar complex 
$$
\Om^sM=M\ox_{E(3)_*}\Sg 3\ox_{E(3)_*}\cdots \ox_{E(3)_*}\Sg 3 \lnr{cobar}
$$
($s$ factors of $\Sg 3$)
with well known differentials $d_r\cln \Om^rM\to \Om^{r+1}M$ (see \kko d).

The cohomology $H^tM^s_n$ of the monochromatic comodules with $s+n=3$ are determined in the following cases
(\cf \cite[6.3.12. Th., 6.3.14. Th.]{r:book}, \cite[Th.~5.10]{mrw}) :
\skr{
H^0M^0_3\eqs K(3)_*\\
H^1M^0_3\eqs K(3)_*\{h_0,h_1,h_2,\ze_3\}\\
H^2M^0_3\eqs K(3)_*\{g_i, k_i, b_i, h_i\ze_3\mid i\in\Z/3\}\qand\\
H^0M^1_2\eqs \ds K(2)_*/k(2)_*\o+ \Op_{i\ge 0, s\in\Zpp}k(2)_*/(v_2^{a_i})\{x_i^s/v_2^{a_i}\},
}{HM}
The cohomology groups
$H^*M^0_3$ and $H^0M^2_1$ are also determined  by Ravenel \cite[6.3.34.~Th.]{r:book} and Nakai \cite n, respectively.
Here, 
$$
k(2)_*=\Z/p[v_2],\ak  K(2)_*=\Z/p[v_2,v_2^{-1}]\qand  K(3)=\Z/p[v_3,v_3^{-1}].
$$
$(K(3)_*=E(3)_*/I_3=M^0_3)$.
The elements
$x_i(=x_{3,i})$ are introduced in \cite[(5.11)]{mrw} such that $x_i\cg v_3^{p^i}$ mod $I_3$ (see Lemma \ref{xi}), and the generators $h_i$, $\ze_3$, $g_i$, $k_i$ and $b_i$ are defined by cocycles in the cobar complex $\Om^*E(3)_*/I_3$ as follows: 
$$
h_i=\hc{\tp i},\ak \ze_3=\hc{Z}, \ak g_i=\hc{G_i}, \ak k_i=\hc{K_i}\qand b_i=\hc{\b1i}.
\lnr{gens}
$$
Hereafter, $\hc x$ denotes the cohomology class represented by a cocycle $x$, and the representatives in \kko{gens} are defined by

\skr{ 
Z\eqs -v_3^{-1}ct_3+v_3^{-p}t_3^p+v_3^{-p^2}t_3\pt2-v_3^{-p}t_1^pt_2\pt2 \\
G_i\eqs t_1\pt i\ox t_2\pt i+\dfrac12t_1^{2p^i}\ox\tp{i+1}\\
K_i\eqs t_2\pt i\ox t_1\pt {i+1} +\dfrac12\tp i\ox t_1^{2p^{i+1}}\qand\\
b_{1,i}\eqs\ds\sum_{k=1}^{p-1}\frac1p\C pk t_1^{kp^i}\ox t_1^{(p-k)p^i}.  
}{coc}
Here, $ct_3$ is the Hopf conjugation of $t_3$ (see Lemma \ref{ct}).
We notice that $G_i$, $K_i$ and $\b1i$ are also cocycles of $\Om^*\B3/I_2$, and of $\Om^*BP_*/I_2$ in \cite[(1.9)]{mrw}.

We introduce integers $e(n)$, $a_n$, $j_{s,n}$ and $j'_{s,n}$ for integers $n\, (\ge 0)$ and $s$ by
\skr{
e(n)\eqs \frac{p^n-1}{p-1}\ak\mbx{for $n\ge 0$,}\\
a_n\eqs \cass{1&\mbox{for $n=0$,}\\ p^n+\frac{p^{n-1}-1}{p+1}&\mbox{for odd $n\ge 1$, and}\\ p^n+p\frac{p^{n-2}-1}{p+1}&\mbox{for even $n \ge 2$}}
}{ea}

\vsp{-.15in}

\skc{
j_{s,n}\eqs \cass{2&\mbox{for $s\in\Z_0$ and $n=0$}\\
2p^2- p+1&\mbox{for $s\in\Z_0$ and $n=2$}\\
2a_{n}+\O 1&\mbox{for $s\in\Z_0$, even $n\ge4$}\\
a_{n+2}-a_{n+1}&\mbox{for $s\in\Z_1$ and even $n\ge 0$ }\\
p+1&\mbox{for $s\in\Zpp$ and $n=1$ }\\
e(3)p^{n-2}-p+1&\mbox{for $s\in\Zpp$ and odd $n\ge 3$ }
}
}1{j0}

\vsp{-.2in}

\skc{
j'_{s,0}\eqs\cass{2&\mbox{for $p\nmid s(s-1)$}\\ 2p&\mbox{for $s=tp+1$ and $p\nmid t(t-1)$}\\
p^2+1 &\mbox{for $s=tp^2+1$ and $p\nmid t$}\\
a_n+p &\mbox{for $s=tp^n+1$ with  $n\ge 2$ and $p\nmid t$}\\
a_n+1 &\mbox{for $s=tp^n+e(n)$ with even $n\ge 2$ and $p\nmid (t-1)$}\\
a_n+2 &\mbox{for $s=tp^n+e(n)$ with odd $n> 2$ and $p\nmid (t-1)$}
}
}2j

\vsp{-.15in}

\skc{
j_{s,n}'\eqs\cass{2p&\mbox{for $s\in\Z_0$ and $n=1$ }\\ 
2pa_{n-1}+p&\mbox{for $s\in\Z_0$ and odd $n\ge 3$ }\\
pa_{n+1}-pa_{n}&\mbox{for $s\in\Z_1$ and odd $n\ge 1$  }\\
p^2+p&\mbox{for $s\in\Zpp$ and $n=2$ }\\
e(3)p^{n-2}-1+\O 1&\mbox{for $s\in\Zpp$ and even $n\ge 4$.}\\
}
}3{j1}

\vsp{-.05in}

\noindent
Here, $\O 1=0$ if $p\ge 5$ and $=1$ if $p=3$, $\Z_i$'s are the subsets of the integers $\Z$ defined in \kko{Zi}, and  the integers $a_n$ are $a_{3,n}$ in \cite[(5.13)]{mrw}. 
We note that 
$$
a_n+a_{n-1}=e(3)p^{n-2}-1\ (n\ge2)\qand p^n+a_{n-2}-p^{n-3}=a_n \ (n\ge3).\lnr{inta}
$$

\thm{H1}{Let $p$ be an odd prime.
$H^1M^1_2$ is the direct sum of $k(2)_*$-module $B_\infty=K(2)_*/k(2)_*\{h_0, h_1,$ $ \wt\ze_2, \ze_3\}$ and $k(2)_*$-cyclic modules generated by

\vsp{-.1in}

\sk{
(\ze_3)_{sp^n/a_n}&\!\mbox{for $(s,n)\in \Zpp\x \N$,}\\ 
\kh0{sp^n/j_{s,n}}&\!\mbox{for $(s,n)\in \LR{(\Z_0\cup \Z_1)\x 2\N}\cup \LR{\Zpp\x \O{2\N}}$,}\\
\kh1{sp^n/j_{s,n}'}&\!\mbox{for $(s,n)\in \LR{(\Z_0\cup \Z_1)\x \O{2\N}}\cup \LR{\LR{\Zpp\x {2\N}}\setminus \{(1,0)\}}$, and}\\
\LE{\hsp{-.7in}\kh2{tp-1/p-1}}&\!\mbox{for $t\in\Z$.}
}
}
We note that there is a little difference between the cases for $p\ge 5$ and $p=3$.
In the theorem, $\wt \ze_2(=(h_1)_1)$ denotes  the homology class of $z$ in \kko{At} (see also \kko{H2M02}),  the generators $(\xi)_{s/j}$ for $\xi=[X]$ in $H^1M^0_3$ denote  

\vsp{-.05in}

\naka{
$
(\xi)_{s/j}=\hc{X/v_2^j+\cdots}
$
}

\vsp{-.05in}

\noindent
for a  cocycle $X/v_2^j+\cdots$ of the cobar comlex $\Om^1M^1_2$ with an element $\cdots$ killed by $v_2^{j-1}$. The element $v_2$ acts on $(\xi)_{s/j}$ by 
$$
v_2(\xi)_{s/j}=(\xi)_{s/j-1} \qand v_2(\xi)_{s/1}=0, \lnr{v2act}
$$
and so, $(\xi)_{s/j}$ generates a cyclic $k(2)_*$-module isomorphic to $k(2)_*/(v_2^j)$:
\naka{$
k(2)_*\{(\xi)_{s/j}\}\cong k(2)_*/(v_2^j).
$}

\section{Proofs of Theorems \ref{H1} and \ref{main1}}

\vsp{-.05in}

In this section, we assume Lemma \ref{key},  which will be verified by a routine calculation in  section six, and prove  Theorems \ref{H1} and \ref{main1}.

\vsp{-.05in}

\subsection{Proof of Theorem \ref{H1}}

For the monochromatic comodules defined in section two, we have a short exact sequence
$$
0\to M^0_3\xar{\et}M^1_2\xar{v_2}M^1_2\to 0, \lnr{ses1}
$$
\noindent
where $\et(x)=x/v_2$ (\cf \cite[(3.10)]{mrw}), which induces
the long exact sequence
%
%
$$
 \cdots\to H^0M^1_2\xar{\de_0}H^1M^0_3\xar{\et_*} H^1M^1_2\xar{v_2} H^1M^1_2
\xar{\de_1}  H^2M^0_3\to \cdots.
\lnr{les}
$$
From  \cite[(5.18)]{mrw}, we read off the following:


\prop{coker}{
  The cokernel of $\de_0\cln H^0M^1_2\to H^1M^0_3$ is a $\Z/p$-module generated by
$(h_0)_0$, $(h_1)_0$, 
\AR{clcl}{
(h_0)_{sp^{2k}}&\qquad  s\in \Z_0\cup \Z_1,
&(h_0)_{tp^{2k+1}}&\qquad t\in \Zpp, \\ 
(h_1)_{tp^{2k}}&\qquad t\in \Zpp, 
&(h_1)_{sp^{2k+1}}&\qquad  s\in \Z_0\cup \Z_1,\\ 
(h_2)_{tp-1}&\qquad t\in \Z,\ak\qand \hsp{-.1in}
&(\ze_3)_t&\qquad t\in \Z
}
for $k\ge 0$.
Here, $\Z_i$ is a subset of $\Z$ given in \kko{Zi}, and $(\xi)_s=v_3^s\xi$ for $\xi\in\{h_i,\ze_3\mid i\in\Z/3\}$.
}


Let $(x)_s\in \Om^1\B3$ denote a cochain satisfying
$$
(x)_s\cg v_3^sx \mod I_3.
$$
\lem{key}{
There exist following cochains in $\Om^1\B3/I_2$: 
\begin{enumerate}
\item $\kt {}{sp^{2k}}$ and $\kt {p}{sp^{2k+1}}$ for $s\in \Z_0$ such that
\sk{
d_1(\kt{}{sp^{2k}})\cgs \cass{s(s+1)v_2^2v_3^{s-1-p}G_2&k=0,\\
s(s+1)v_2^{2p^2-p+1}v_3^{sp^2-2p}G_1&k=1,\\
-3s(s+1)v_2^{2a_{2k}}v_3^{(sp-2)p^{2k-1}}K_0 
&k\ge 2,\ p\ge 5,\\
-2s(s+1)v_2^{2a_{2k}+1}v_3^{3^{2k-1}(3s-2)}(\b10+t_1^p\ox Z')&k\ge 2,\ p=3; and}\\
d_1(\kt{p}{sp^{2k+1}})\cgs \cass{s(s+1)v_2^{2p}v_3^{sp-2}G_0&k=0,\\
s(s+1)v_2^{2pa_{2k}+p}v_3^{(sp-2)p^{2k}}\b11 &k\ge 1.
}
}
\item $\kt {}{sp^{2k}}$ and $\kt {p}{sp^{2k+1}}$ for $s=tp^2-1\in \Z_1$ such that
\sk{
d_1(\kt{}{sp^{2k}})\cgs v_2^{a_{2k+2}-a_{2k+1}}v_3^{(tp-1)p^{2k+1}}\b10 \qand\\
d_1(\kt{p}{sp^{2k+1}})\cgs v_2^{pa_{2k+2}-pa_{2k+1}}v_3^{(tp-1)p^{2k+2}}\b11\mbx{for $k\ge 0$.}
}
\item $\kt {}{sp^{2k+1}}$ and $\kt {p}{sp^{2k}}$ for $s\in \Zpp$ such that
\sk{
d_1(\kt p{tp^k+1})\cgs \cass{t(t-1)v_2^{2p}v_3^{tp-1}G_0&k=1,\\
-tv_2^{p^2+1}v_3^{(tp-1)p}G_1&k=2,\\
-2tv_2^{a_k+p}v_3^{(tp-1)p^{k-1}}G_0&\text{odd $k\ge 3$,}\\
2tv_2^{a_k+p}v_3^{(tp-1)p^{k-1}}K_0&\text{even $k\ge 4$;}
}\\
d_1(\kt p{tp^k+e(k)})\cgs\cass{-(t-1)v_2^{a_k+1}v_3^{tp^k+pe(k-2)}G_1&\text{even $k\ge 2$,}\\
-(t-1)v_2^{a_k+2}v_3^{tp^k+pe(k-2)}b_{1,1}&\text{odd $k\ge 3$;}
}\\
d_1(\kt{p}{sp^{2k}})\cgs \cass{s(s-1)v_2^2v_3^{s-2}K_1&k=0,\\
-sv_2^{p^2+p}v_3^{sp^2-p-1}K_0&k=1,\\
-3sv_2^{e(3)p^{2k-2}-1}v_3^{(sp^2-p-1)p^{2k-2}}K_0 
&p\ge 5,\ k\ge 2,\\
-sv_2^{3^{2k-2}e(3)}v_3^{(9s-4)3^{2k-2}}(\b10+Z'\ox t_1^p)
&p=3,\ k\ge 2; \text{ and}}\\
d_1(\kt{}{sp^{2k+1}})\cgs \cass{-sv_2^{p+1}v_3^{(s-2)p}K_2&k=0,\\
sv_2^{e(3)p^{2k-1}-p+1}v_3^{(sp^2-p-1)p^{2k-1}}\b11&k\ge 1.}
}

\item $\kt {p^2}{tp-1}$  such that $d_1(\kt {p^2}{tp-1})\cg v_2^{p-1}v_3^{tp-p}b_{1,2}$.
\end{enumerate}
Here,  $G_i$, $K_i$ and $\b1i$ are the cocycles of $\Om^2\B3/I_2$ 
in \kko{coc}, 
$Z'$ is an element in Lemma \ref{xi}, and $x\cg v_2^ay$ denotes the congruence modulo $J_{a+1}$. 
}


Let $d_1((x)_t)\cg v_2^j y$ mod $J_{j+1}$ be a congrucence in Lemma \ref{key}. 
Then, $\de_1((\hc x)_{t/j})=\hc y$
for  the connecting homomorphism $\de_1$ in \kko{les}.
Here, $(\hc x)_{t/j}(=\hc{(x)_t/v_2^j})\in H^1M^1_2$ denotes the cohomology class of the cocycle 
$(x)_t/v_2^j$ of $\Om^1M^1_2$.
Thus, the cochains in Lemma \ref{key} give rise to elements $\kh0{sp^r/j_{s,r}}$ and $\kh1{sp^r/j'_{s,r}}$ of $H^1M^1_2$ as well as the $\de_1$-images of them.
Furthermore, we have elements
$$
(\ze_3)_{tp^n/a_n}=x_n^t\ze_3/v_2^{a_n}\in H^1M^1_2
$$
for the elements $x_n(=x_{3,n})$ introduced in \cite[(5.11)]{mrw} (see Lemma \ref{xi})  with
 $$
 \de_1((\ze_3)_{tp^n/a_n})=\cass{(h_0\ze_3)_{(tp-1)p^{n-1}}&\text{ $n$ is odd,}\\
 (h_1\ze_3)_{(tp-1)p^{n-1}}&\text{ $n$ is even}} \lnr{dhz3}
$$
 by \cite[(5.18)]{mrw} (or Lemma \ref{xi}).
We notice that as a $k(2)_*$-module,
 $K(2)_*/k(2)_*\{\xi\}=\Z/p\{(\xi)_{0/j}\mid j\ge 1\}$ with $v_2(\xi)_{0/j}=(\xi)_{0/j-1}$ and $v_2(\xi)_{0/1}=0$ (see \kko{v2act}). 

Let $B$ be the $k(2)_*$-module of the theorem.
Each direct summand of $B$ is a submodule of $H^1M^1_2$, which defines a $k(2)_*$-module map $f\cln B\to H^1M^1_2$.
Furthermore, assigning $(\xi)_{s/1}\in B$ to the generator $(\xi)_s$ of the cokernel of $\de_0$, we have a homomorphism $\O\et_*\cln H^1M^0_3\to B$ by Proposition \ref{coker}.
These homomorphisms fit in the commutative diagram 
\xyR{10}{
H^0M^1_2\ar[r]^-{\de_0}\ar@{=}[d]&H^1M^0_3\ar[r]^-{\O\et_*}\ar@{=}[d]&B\ar[r]^-{v_2}\ar[d]^-f&B\ar[r]^-{\de'_1}\ar[d]^-f&H^2M^0_3\ar@{=}[d]\\
H^0M^1_2\ar[r]^-{\de_0}&H^1M^0_3\ar[r]^-{\et_*}&H^1M^1_2\ar[r]^-{v_2}&H^1M^1_2\ar[r]^-{\de_1}&H^2M^0_3
,}
where we define $\de'_1$ by $\de_1f$.
It suffices to show that the upper sequence is exact by \cite[Remark 3.11]{mrw}.
By the definition of $B$, the subsequence $H^0M^1_2\xar{\de_0}H^1M^0_3\xar{\O\et_*}B\xar{v_2}B$ is exact 
and the composite $B\xar{v_2}B\xar{\de'_1}H^2M^0_3$ is zero.

Suppose that the $\de'_1$-images of the generators are linearly independent, and
take $\xi\in \Ker \de'_1$ to be a homogeneous element.
Then,
\AR c{
\xi=\sum_k c_k\xi_k \mbx{for  generators $\xi_k$ of $B$ and scalars $c_k\in k(2)_*$, and} \\
0=\de'_1(\xi)=\sum_k \O c_k\de'_1(\xi_k)
}
for the image $\O c_k$ of $c_k$ under the projection $k(2)_*\to \Z/p$ sending $v_2$ to zero.
Since  $\de'_1(\xi_k)$'s are linearly independent, we see $\O c_k=0$, and
so we have $c_k'\in k(2)_*$ such that $c_k=v_2c_k'$. 
Therefore,
$$\textstyle
\xi=\sum_kv_2c'_k\xi_k\in \im v_2,
$$
and we see the upper sequence of the above diagram is exact if the $\de'_1$-images of the generators are linearly independent.

The $\de'_1$-image is a $\Z/p$-submodule of $H^2M^0_3$ in \kko{HM} generated by the generators of the form $(\rh)_{s}$ for $\rh\in\{g_i,k_i,b_i, h_i\ze_3\mid i\in\Z/3\}$ by Lemma \ref{key} and \kko{dhz3}. 
Moreover,
Lemma \ref{key} and \kko{dhz3} show that the $\de'_1$-image of each generator $\xi_k$ has the only one pattern of form $(\rh)_s$ 
 except for 
%
 \AR{|c||c|c|c|c|}{
\cline{1-3}
g_0&(g_{0})_{sp-2}
&(g_0)_{(tp-1)p^{2n}}\\
\cline{1-4}
g_1&(g_{1})_{(sp-2)p}& (g_1)_{(tp-1)p}&(g_1)_{tp^{2n}+pe(2n-2)} \\
\cline{1-4}
k_0&(k_0)_{(sp-2)p^{2n-1}}&(k_0)_{(tp-1)p^{2n-1}}&(k_0)_{(sp^2-p-1)p^{2n}}\\
\hline
b_1&(b_1)_{(sp-2)p^{2n}}&(b_1)_{(tp-1)p^{2n+2}}&(b_1)_{t{p^{2n+1}}+pe(2n-1)}&(b_1)_{(sp^2-p-1){p^{2n-1}}}\\
\hline\hline
b_0&(p=3)&(b_0)_{3^{2n-1}(3s-2)}&(b_0)_{3^{2n+1}(3t-1)}&(b_0)_{3^{2n-2}(9s-4)}\\
\hline
}
These results show the independency of the $\de'_1$-images as desired.
\q

\subsection{Proof of Theorem \ref{main1}}
Let $\de^0_2\cln H^*N_2^{1}\to H^{*+1}N^0_2$ be the connecting homomorphism associated to the short exact sequence \kko{ses:ch}, and
consider the diagram
\xyR{10}{
&H^2M^0_2\ar[r]^-{(\ka^0_2)_*}&H^2N^1_2\ar[r]^-{\de^0_2}\ar[d]^-{\io^1_2}&H^3N^0_2=E_2^{3}(V(1))\\
H^1M^1_2\ar[r]^-{\de_1}&H^2M^0_3\ar[r]^-{\et_*}&H^2M^1_2
}
of exact sequences
for $\de_1$ in \kko{les}.
The connecting homomorphism $\O\de_j$ associated to \kko{ses2} 
is factorized into the composite $\O\de_j\cln H^s\A/J_j\xar{\wh \io_j} H^sN^1_2\xar{\de^0_2}H^{s+1}N^0_2$
for the homomorphism $\wh\io_j$ given by $\wh\io_j(x)=x/v_2^j$.
It follows that
$$
 \O \ga''_{sp^r/j}=\de^0_2(v_3^{sp^r}/v_2^j)\in H^1N^0_2=E_2^{1}(V(1))\mbx{for $v_3^{sp^r}/v_2^j\in H^0N^1_2$.} \lnr{ga}
$$
Since $\de^0_2$ is a $k(2)_*$-module map, we have 
$$
v_2^{j-1}\O \ga''_{sp^r/j}=v_2^{j-1}\de^0_2(v_3^{sp^r}/v_2^j)=\de^0_2(v_2^{j-1}v_3^{sp^r}/v_2^j)
=\de^0_2(v_3^{sp^r}/v_2)=\O \ga''_{sp^r}. \lnr{v2ga}
$$

It is well known that
 $$
 \O\be_1=-b_0=\left[-\b10\right], \qand  \O\be_2=2k_0=\left[2K_0\right]\in H^2N^0_3
 $$
 for the cocycles $\b10$ and $K_0$ in \kko{coc} 
(\cf \cite[Lemma 4.4]{os}).
This defines elements $v_3^{sp^r}\O\be_i/v_2\in H^2N^1_2$ for $i=1,2$, and 
$$
\de^0_2(v_3^{sp^r}\O\be_i/v_2)\usk{ga}=\ga''_{sp^r}\O\be_i\in E_2^3(V(1)). 
$$
We also see that for $v_3^{sp^r}\O\be_i\in H^2M^0_3$,
$$
\et_*(v_3^{sp^r}\O\be_i)=\io^1_2 (v_3^{sp^r}\O\be_i/v_2)\in H^2M^1_2.
$$
From 
 Lemma \ref{key}, we read off that the elements $v_3^{sp^r}\O\be_1=-(b_0)_{sp^r}$ and $v_3^{sp^r}\O\be_2=2(k_0)_{sp^r}\in H^2M^0_3$ have a possibility to be in the image of
 $\de_1$ if
\AR l{
\mbox{(a)\ak $p\ge 5$ and  $(s,r)\in \LR{\Z_1^+\cup \Z_2^+}\x \O{2\N}$, or}\\ 
\mbox{(b)\ak $p=3$ and $(s,r)\in   \LR{\O{\N_{1}}\x 2\N_{>0}}\cup \LR{\LR{\Z_1^+\cup \Z_2^+}\x \O{2\N}}\cup \LR{\O{\N_2}\x\O{2\N}_{>1}}$, } 
}
and if 
\AR l{
\mbox{(a)\ak $p\ge 5$ and $(s,r)\in   \LR{\O{\N_{1}}\x 2\N}\cup \LR{ \Z_2^+\x \O{2\N}_{>1}}\cup \LR{\O{\N_2}\x\O{2\N}_{>1}}$, or} \\
\mbox{(b)\ak $p=3$ and  $(s,r)\in  \LR{\O{\N_{1}}\x \{0\}}\cup \LR{\Z_2^+\x \O{2\N}_{>1}}$,}\\ 
}
respectively. 
Here, $\O{\N_i}=\Nt p\setminus \N_i$ for $i=1,2$. 
Therefore, if a pair $(s, r)$ satisfies the condition of the theorem, then the element $v_3^{sp^r}\O\be_i$ is not in the image of $\de_1$,
and survives to $\io^1_2(v_3^{sp^r}\O\be_i/v_2)$ under the homomorphism $\et_*$.
Thus, $v_3^{sp^r}\O\be_i/v_2\ne 0\in H^2N^1_2$ under the conditions.

Ravenel determined in \cite[6.3.24.~Th.]{r:book} and \cite[(3.2)~Th.]{r:coh} that
$$
H^2M^0_2=\cass{K(2)_*\{h_0\wt\ze_2, h_1\wt\ze_2, b_0, b_1, \xi\}&p=3\\ K(2)_*\{ h_0\wt\ze_2, h_1\wt\ze_2, g_0, g_1\} &p\ge 5},
\lnr{H2M02}
$$
where $\wt\ze_2=v_2^{p+1}\ze_2=\hc{-z}$ for $\ze_2$ in \cite[Prop. 3.18)]{mrw} and $z$ in \kko{At}.
This shows that the elements $v_3^{sp^r}\O\be_i/v_2$ for $i=1,2$ are not in the image of $(\ka_2^0)_*$,
and hence  survive to $\ga''_{sp^r}\O\be_i\in E_2^3(V(1))$.
Moreover, $\ga''_{sp^r/j}\O\be_i\ne 0\in E_2^3(V(1))$ if $v_2^{j-1}\ga''_{sp^r/j}\O\be_i\usk{v2ga}=\ga''_{sp^r}\O\be_i$ is not zero.
\q

\section{Some cochains in the cobar complex $\Om^*E(3)_*$}

In the rest of this paper, we consider $\Sg3$-comodules whose structure maps are induced from the right unit map $\eR\cln E(3)_*\to \Sg3$.
We consider the cobar complex $\Om^*M$ of a comodule $M$ in \kko{cobar}, whose differentials  are given by
\skr{
d_0(v)\eqs \eR(v)-v\in \Om^1E(3)_*, \qand\\
d_1(x)\eqs 1\ox x-\De(x) +x\ox 1\in \Om^2E(3)_*\\ 
}d
 for $v\in \Om^0E(3)_*=E(3)_*$ and $x\in \Om^1E(3)_*=\Sg 3 $.
 For the differentials $d_0$ and $d_1$, we have  relations (\cf \cite[(2.3.2)]{st}):
 \skr{
d_0(vv')\eqs vd_0(v')+d_0(v)\eR(v'),\\
d_1(vx)\eqs d_0(v)\ox x+vd_1(x) \\
d_1(xy)\eqs -x\ox y-y\ox x+d_1(x)\De y+(x\ox 1+1\ox x)d_1(y) \qand\\
d_1(x\eR(v))\eqs d_1(x)(1\ox \eR (v))-x\ox d_0(v) 
}{deriv}
for $v,v'\in E(3)_*$ and $x,y\in \Sg3$.
A formula for the Hopf conjugation $c\cln \Ga\to \Ga$ is given in \cite[(3)]{r:BP}, 
and implies immediately the following:

\lem{ct}{
The Hopf conjuataion $c\cln E(3)_*(E(3))\to E(3)_*(E(3))$ acts as\\[-4mm]
\AR c{ct_1=-t_1,\ak
ct_2=t_1^{p+1}-t_2, \qand ct_3\cg t_2t_1\pt2-t_1ct_2^p-t_3 \mod I_2.
}
}

For the right unit $\eR\cln \A\to \Ga$, 
we have a well known formula\\[-4mm]
$$
\eR(v_n)\cg v_n+\vm \tp{n-1}-\vm^pt_1\mod I_{n-1}. \leqno{(\nr)\label{L}(\mbox{\cite[(11)]{r:BP}})  }
$$
A routine calculation using \kko d and \kko{L} shows the following: 
\lem{si}{
Put $\si_n=\sum_{k=0}^{n-1}v_2^{p^{2k}a_{2n-2k-1}-p^{2k+1}}v_3^{p^{2k}}\in \B3$.
Then,
$$
d_0(\si_n)\cg v_2^{p^{2n-2}}\tp{2n}-v_2^{a_{2n-1}}t_1\mod I_2.
$$
}

In $\Sg3$, $\eR(v_4)=0=\eR(v_5)$, which give rise to relations\\[-4mm]
\skr{
v_3\tp3\cgs t_1\eR(v_3)^p-v_2t_2\pt2+v_2\pt2t_2\qand \\
v_3t_2\pt3\cgs t_2\eR(v_3)\pt2-v_2t_3\pt2-v_2w^p+v_2\pt3t_3\mod I_2
}{t3=t} 

\vsp{-1mm}
\noindent (\cf \cite[(12),  (16)]{r:BP}, \cite[4.3.21.~Cor.]{r:book}), where
$w\in \Sg3$ ($=w_1(v_3,v_2\tp2,-v_2^pt_1)$ in \cite[4.3.21.~Cor.]{r:book}) is an element defined by\\[-4mm]
$$
pw=v_3^p+v_2^p\tp3-v_2\pt2t_1^p+y^p-\eR(v_3)^p\lnr{w13}
$$

\vsp{-1mm}
\noindent
for $y\in (p,v_1)$ in $\eR(v_3)=v_3+v_2\tp2-v_2^pt_1+y$ (see \kko{L}).

The diagonal $\De\cln E(3)_*(E(3))\to E(3)_*(E(3))\ox_{E(3)_*}E(3)_*(E(3))$ of the Hopf algebroid $\Sg 3$ acts on the elements $t_i$ and $ct_i$ as follows:\\[-4mm]
\skr{
\De(t_1)\eqs t_1\ox 1+1\ox t_1,\\
\De(t_2)\cgs t_2\ox 1+t_1\ox t_1^p+1\ox t_2-v_1\b10\mod (p,v_1^2), \qand\\
\De(t_3)\cgs t_3\ox 1+t_2\ox t_1^{p^2}+t_1\ox t_2^p+1\ox t_3-v_2b_{1,1}\mod I_2\\
\De(t_4)\cgs t_4\ox 1+t_3\ox t_1^{p^3}+t_2\ox t_2^{p^2}+t_1\ox t_3^p+1\ox t_4-v_3b_{1,2}\mod I_3
}{Det}

\vsp{-.5mm}
\noindent
(\cf \cite[Th.~8]{r:BP}, \cite[4.3.15. Cor.]{r:book}), and so\\[-4mm]
\skr{
d_1(ct_2)\cgs -t_1^p\ox t_1,\\
d_1(ct_3)\cgs ct_2^p\ox t_1+\tp2\ox ct_2-v_2b_{1,1}\mod I_2\qand\\
d_1(ct_4)\cgs t_1^{p^3}\ox ct_3-ct_2^{p^2}\ox ct_2+ct_3^p\ox t_1-v_3b_{1,2}\mod I_3,
}{Dec}

\vsp{-.5mm}
\noindent
since $\De(cx)=(c\ox c)T\De(x)$ for the switching map $T$ given by $T(x\ox y)=y\ox x$,
 where  $\b1k$ is the cocycle in \kko{coc}.

The fact 
$
d_1(\tp{k+1})\cg -p\b1k 
$ mod $(p^2)$ implies not only that
the cochain $\b1k\in \Om^2E(3)_*/(p)$ is a cocycle, but also the following lemma.

\lem w{
The cochain $w$ in \kko{w13} satisfies
$$
w\cg -v_2v_3^{p-1}t_1\pt2 \mod J_{2}\qand
d_1(w)\cg  -v_2^{p}b_{1,2} + v_2^{p^2}\b10 \mod I_2.
$$
}

\cor b{
Put $W_n=\sum_{i=0}^{n-1} v_2^{p^{2i}a_{2n-2i}-p^{2i+2}} w\pt {2i}$.  
Then,
$$
d_1(W_n)\cg -v_2\pt{2n-1} \b1{2n}+v_2^{a_{2n}}\b10 \mod I_2.
$$
}

We generalize the relations \kko{t3=t} and obtain the following proposition
from \cite[(4.3.1), 4.3.11 Lemma]{r:book} and \cite[Th.~1]{r:BP} (\cf \cite[Prop.~2.1]{seR}): 

\prop{3.5}{
There exist elements $T_{n}$ for $n\ge 0$ satisfying
$T_{n}\cg t_n^p$ {\rm mod} $I_3$ and   
$$
v_2\pt{k+1}t_{k+1}+t_{k}\eR(v_3)\pt{k}\cg v_1T_{k+2}+v_2T_{k+1}^p+v_3T_{k}\pt2\mod (p, v_1^2)
$$
for $k\ge 0$,
In particular, $T_0=1$, $T_1\cg t_1^p$, $T_2\cg t_2^p$ and $T_3\cg t_3^p+w$ {\rm mod} $I_2$.
}

\p
We begin with recalling some notations from \cite[\S4.3]{r:book}.
For a sequence $J=(j_1,j_2,\dots, j_m)$ of positive integers, we set $|J|=m$ and $\|J\|=\sum_{i=1}^mj_i$, and  an element $v_J\in\B3$ is defined recursively by
$v_{(j,J)}=v_{j}v_J\pt j$.
Let
$w_k(S)$ for a set $S$ be symmetric polynomials of degree $p^n$ such that
$w_0(S)=\sum_{x\in S}x$ and $\sum_{x\in S}x\pt n=\sum_{k=0}^np^kw_k(S)\pt{n-k}$.
We then define sets $S_n$ out of a set $S=\{a_{i,j}\}$ recursively by
$$
S_n= \{a_{i,j}\mid i+j=n\}\cup\bigcup_{|J|>0}\{v_Jw_{|J|}(S_{n-\|J\|})\pt{\|J\|-|J|}\}.
$$

\vsp{-.1in}

\noindent
By \cite[(4.3.1), 4.3.11 Lemma]{r:book}, we see

\vsp{-.13in}

$$
w_0(C_n)\cg {\sum_{i+j=n}}^{\!\!\!F}t_i\eR(v_j)\pt i\cg {\sum_{i+j=n}}^{\!\!\!F}v_it_j\pt i\cg w_0(D_n) \mod (p)\leqno{(\nr)\label{wC=wD}}
$$

\vsp{-.1in}

\noindent
for the sets
$$
C=\{t_i\eR(v_{j})\pt i\} \qand D=\{v_it_{j}\pt i\}
$$
In $\Sg3$, put
$$
w(S_n)=\sum_{J}v_J^pw_{|J|+1}(S_{n-\|J\|})\pt{\|J\|-|J|} 
\qand 
T_n=t_n^p-w(C_n)+w(D_n).
$$
Then, the proposition follows from \kko{wC=wD} and the congruences
\sk{
w_0(C_n)\cgs v_2\pt{n-2}t_{n-2}+t_{n-3}\eR(v_3)\pt{n-3}+v_1w(C_{n-1})+v_2w(C_{n-2})^p+v_3w(C_{n-3})\pt2\\
w_0(D_n)\cgs v_1t_{n-1}^p+v_2t_{n-2}\pt2+v_3t_{n-3}\pt{3}+v_1w(D_{n-1})+v_2w(D_{n-2})^p+v_3w(D_{n-3})\pt2\\
}
seen by the relation 
$$
v_{(k,J)}w_{|(k,J)|}(S_{n-\|(k,J)\|})\pt{\|(k,J)\|-|(k,J)|}=v_kv_J\pt kw_{|J|+1}(S_{n-k-\|J\|})\pt{\|J\|-|J|+k-1}.
$$

\vsp{-.2in}
\q

\lem{3.6}{For $n\ge 0$,

\vsp{-.15in}

$$
\eR(v_2^{p-1}v_3^{e(n)})\cg \sum_{i=0}^n(-1)^{n-i}v_2^{p^{i+1}e(n-i)+p-1}v_3^{e(i)}t_{n-i}\pt i-v_2^pw_n^p+v_1v_2^{p-2}w_{n+1} \mod (p,v_1^2).
$$
Here,

\vsp{-.2in}

$$
w_n=\sum_{i=1}^n(-1)^iv_2^{e(i-1)}T_i\eR(v_3^{p^{i-1}e(n-i)}). \lnr{wn}
$$
}

\p In this proof, every congruence is considered modulo $(p,v_2^2)$.
By Proposition \ref{3.5}, we have $t_k\eR(v_3\pt k)\cg \wt T_k-v_2\pt{k+1}t_{k+1}$ for $\wt T_k=v_1T_{k+2}+v_2T_{k+1}^p+v_3T_{k}\pt 2$, which implies inductively

\vsp{-.2in}

$$
t_1\eR(v_3^{pe(n)})\cg -\sum_{i=1}^n (-1)^iv_2^{p^2e(i-1)}\wt T_i\eR(v_3^{p^{i+1}e(n-i)})+(-1)^nv_2^{p^2e(n)}t_{n+1},
$$

\vsp{-.1in}

\noindent
and hence

\vsp{-.2in}

\skr{
t_1\eR(v_3^{pe(n)})\cgs -v_1v_2^{-p-1}w_{n+2}+v_2^{1-p}w_{n+1}^p-v_3w_{n}\pt2+(-1)^nv_2^{p^2e(n)}t_{n+1}\\
&\hsp{-.3in}-v_1v_2^{-p-1}(t_1^p\eR(v_3)-v_2t_2^p)\eR(v_3^{pe(n)})+v_2^{1-p}\tp2\eR(v_3^{pe(n)}).
}{t1R} 
\indent
Now we prove the lemma by induction. For $n=0$, it follows from the facts: $\eR(v_2)\cg v_2+v_1t_1^p$  by \kko L and $w_1=-t_1^p$.

Assuming the case for $n$, we obtain the case for $n+1$ from \kko{t1R} and
\sk{
\eR(v_2^{p-1}v_3^{e(n+1)})\cgs v_2^{-p^2+2p-1}v_3\eR(v_2^{p-1}v_3^{e(n)})^p+v_2^{p-1}(v_2\tp2+v_1t_2^p)\eR(v_3^{pe(n)})\\
&-v_2^{2p-1}t_1\eR(v_3^{pe(n)})-v_1v_2^{p-2}t_1^p\eR(v_3^{e(n+1)}),
}
given by $\eR(v_2^{p-1}v_3)\cg v_2^{p-1}(v_3+v_2\tp2-v_2^pt_1+v_1t_2^p)-v_1v_2^{p-2}t_1^p\eR(v_3)$.
Here, $\eR(v_3)$ is given in \cite[(5.7)]{mrw}.
\q

Send the congruence in Lemma \ref{3.6} under $d_1$, and compare the $v_1$-multiples.
Then, we deduce the following corollary (\cf \cite[Prop.~2.3]{seR}).
Indeed, if $v_1v_2^{p-2}d_1(w_{n+1})$ $\cg A+v_1B$ mod $(p,v_1^2)$ for some $A$, $B$ involving no $v_1$, then $A\cg 0$ mod $(p,v_1^2)$ and $v_2^{p-2}d_1(w_{n+1})\cg B$ mod $I_2$. 

\cor{3.7}{For the elements $w_n$ in \kko{wn},

\vsp{-.15in}

$$
d_1(w_{n+1})\cg -\sum_{i=0}^{n-1}(-1)^{n-i}v_2^{p^{i+1}e(n-i)}w_{i+1}\ox t_{n-i}\pt i -(-1)^nv_2^{e(n+1)}\db_n \mod I_2.
$$
Here, $\db_n$ is an element in $d_1(t_n)\cg \da_n+v_1\db_n$ {\rm mod} $(p,v_1^2)$ for  $\da_n$ and $\db_n$ involving no $v_1$.
In particular, $\db_2=\b10$ by \kko{Det}.
}

\vsp{-.05in}

We have the cocycle $z$ in $\Om^1E(3)_*/I_2$:
$$
z=v_3t_1^p+v_2ct_2^p-v_2^pt_2=t_1^p\eR(v_3)-v_2t_2^p+v_2^pct_2=-w_2+v_2^pct_2,
\lnr{At}
$$
which represents the element $ -v_2^{p+1}\ze_2\in H^1M^0_2$ (\cf \cite[Prop.~3.18 c)]{mrw}, \kko{H2M02}). 
In particular,
$$
t_1^p\eR(v_3)\cg z+v_2t_2^p-v_2^pct_2 \mod I_2\lnr{t1eRv3}
$$
We further  have cocycles  $G'_i$ and $K'_i\in \Om^2E(3)_*/I_2$  for $i\in \{0,1,2\}$ defined by
$$
G'_i=ct_2\pt i\ox t_1\pt i+\frac12\tp{i+1}\ox t_1^{2p^i}\qand K'_i=t_1\pt {i+1}\ox ct_2\pt i+\frac12t_1^{2p^{i+1}}\ox\tp i,
\lnr{GK}
$$
which are homologous to $G_i$ and $K_i$ in \kko{coc}, respectively.
Indeed,
$$
d_1(\dg_ i)\cg G_i'-G_i\qand d_1(\dk_ i)\cg K'_i-K_i \mod I_2, \lnr{homo}  
$$
for $i\in \{0,1,2\}$, and for $\dg_i$ and $\dk_i\in \Om^1E(3)_*$ given by
$$
\dg_i=t_1\pt it_2\pt i-\frac12t_1^{p^{i+1}+2p^i}\qand \dk_i=t_1\pt {i+1}t_2\pt i-\frac12t_1^{2p^{i+1}+p^i}. \lnr{dgk}
$$
We also have similar relation
$$
d_1(t_1^pt_2)\cg -(t_1^p\ox t_2+ct_2\ox t_1^p)-2K_0\mod I_2. \lnr{2K0}
$$

\lem{xz}{
In $\Om^1E(3)_*$, put
\AR c{
\om_1=\eR(v_3)t_2-v_2t_3+v_2^pt_1t_2, \ak \om_2=\frac12\eR(v_3)t_1^{2p}-v_2^p\dk_0,\qand\\
\wt\om_2=-w_3-v_2^{pe(2)}t_1^pt_2.
}
Then, modulo $I_2$, 
\sk{
d_1(\om_1)\cgs -t_1\ox z-v_2^2\b11-2v_2^pG_0, \\ 
d_1(\om_2)\cgs -t_1^p\ox z-v_2G_1 +v_2^pK_0,  \qand\\
d_1(\wt\om_2)\cgs v_2\pt2 z\ox t_1^p +2v_2^{p^2+p}K_0+v_2^{e(3)}\b10.
}
}

\p  In this proof, we consider congruences modulo $I_2$.
A routine calculation shows the congruence for $d_1(\om_1)$: 
\sk{
d_1(\eR(v_3)t_2)
\usscgs{\kko{deriv}\\ \kko{t1eRv3}}\!\! -t_1\ox (z+\U{v_2t_2^p}a+\uwave{v_2^pt_2}-\U{v_2^pt_1^{p+1}}c)-t_2\ox (\U{v_2\tp2}b-\U{v_2^pt_1}d)\\[-2mm]
d_1(-v_2t_3)\!\!\ukcgs{Det}\!\! v_2(\U{t_1\ox t_2^p}a+\U{t_2\ox \tp2}b-v_2\b11)\\[-2mm]
d_1(v_2^pt_1t_2)\usscgs{\kko{Det}\\\kko{deriv}} -v_2^p(\uwave{t_1\ox t_2}+\U{t_2\ox t_1}d+\uwave{t_1^2\ox t_1^{p}}+\U{t_1\ox t_1^{p+1}}c),\\[-2mm]
}
in which the underlined terms with the same subscript cancel each other and the wavy underlined terms make $-2v_2^pG_0$. 

For $d_1(\om_2)$, we calculate
$$\textstyle
d_1(\frac12\eR(v_3)t_1^{2p})\!\!\uss{\kko{deriv}\\\kko{t1eRv3}}\cg \!\!- t_1^p\ox (z+\UW{v_2t_2^p}_G\!\!\!\!-\UW{v_2^pct_2}_{K'}\!\!\!\!)-\UW{\frac12 v_2t_1^{2p} \ox \tp2}_G\!\!\!\!\!+\UW{\frac12 v_2^pt_1^{2p}\ox t_1.}_{K'}
$$

\vsp{-1.5mm}
\noindent
Add $d_1(-v_2^p\dk_0)$, and we obtain the desired conguence by \kko{homo}. 

We verify $d_1(\wt\om_2)$ by 
\sk{
d_1(w_3)\urcgs{3.7} -v_2^{pe(2)}w_1\ox t_2+v_2^{p^2}w_2\ox t_1^p-v_2^{e(3)}\b10 \\[-1.5mm]
\usscgs{\kko{At}\\\kko{wn}} -\U{v_2^{pe(2)}(-t_1^p)\ox t_2}a+v_2^{p^2}(-z+\U{v_2^pct_2}b)\ox t_1^p-v_2^{e(3)}\b10\\[-1.5mm]
d_1(v_2^{pe(2)}t_1^pt_2)\ukcgs{2K0} -v_2^{pe(2)}((\U{t_1^p\ox t_2}a+\U{ct_2\ox t_1^p}b)+2K_0).\\
}

\vsp{-.3in}
\q

\section{The elements $x_i$ and deriving elements $y_i$ and $y_i'$}

In \cite[(5.11)]{mrw}, Miller, Ravenel and Wilson introduced elements $x_{3,i}\in v_3^{-1}BP_*$.
We refine them, and define the elements $x_i\in E(3)_*$ by
\AR{l}{
x_i= v_3\pt i\mbx{for $i=0,1,2$,}\hsp{.5in}
x_3= x_2^p-v_2^{p^3-1}v_3^{(p-1)p^2+1},\\
x_4= x_3^p-v_2^{e(2)p^3-p-1}v_3^{(p^2-e(2))p^2+p+1}, \\
x_{2k+1}=x_{2k}^p-v_2^{pa_{2k}-1}x_{2k-1}^{(p-1)p}v_3-v_2^{e(3)p^{2k-1}-e(3)}v_3^{(p^2-e(2))p^{2k-1}+p+1}, \qand\\
x_{2k+2}=x_{2k+1}^p-2v_2^{e(3)p^{2k}-e(3)}v_3^{(p^2-e(2))p^{2k}+p+1}
}
for $k\ge 2$.
\lemc{xi}{\cf \cite[Prop.~3.1]{seR}}{
In $\Om^1E(3)_*$, we have
\sk{
d_0(x_0)\cgs v_2\tp2-v_2^pt_1\mod I_2,\\
d_0(x_1)\cgs v_2^pv_3^{p-1}t_1-v_2^{p+1}v_3^{-1}t_2\pt2\mod J_{2p}, \qand\\
d_0(x_i)\cgs v_2^{a_i}(x_{i-1}^{p-1}\tp{\ep_i}+B_i) \mod J_{e(3)p^{i-2}}\mbx{for $i\ge 2$.}
}
Here,  $\ep_i=\frac{1+(-1)^i}2$, and $B_i$ are as follows
\AR{|c||c|c|c|}{
\hline
i&2&3&2k\\
\hline
B_i&-v_2^pv_3^{c(2)}t_2&v_2^{p^2-p}v_3^{c(3)}(z-v_2^pt_1^{p+1})&v_2^{a_{2k-1}-p}v_3^{c(2k)}(z-v_2^pt_2)\\
\hline
}
\hsp{.38in}
\aR{|c||c|}{
\hline
i&2k+1\\
\hline
B_i&v_2^{a_{2k}-p}v_3^{c(2k+1)}(2z-v_2^pct_2)\\
\hline
}\\

\noindent
for $c(k)=(p^2-p-1)p^{k-2}$.
For $i\ge 4$, add $v_2^{a_{i-1}+1}v_3^{c(i)}Z'$ to $B_{i}$ if we consider the congrucence modulo $J_{e(3)p^{i-2}+1}$.
Here, $Z'$ is a cocycle homologous to $aZ$ for some $a\in\Z/p$.}

\p This follows from a routine calculation:
For $i\le2$, it follows from \kko L and from \kko{t3=t}.

We obtain $d_0(x_3)$ from \kko{t1eRv3}
 and $d_0(v_3^{(p-1)p^2+1})\cg v_3^{(p-1)p^2}(v_2\tp2-v_2^pt_1)-v_2^{a_2}v_3^{(p-1)p^2-p}(t_1^p\eR(v_3)-v_2^pt_2)$ mod $J_{e(3)}$ by \kko{deriv}, \kko L and the congruence on $d_0(x_2)$.
We note that $\eR(v_3^{p+1})=v_3^{p+1}+v_2z^p-v_2\pt2z$ by \cite[(3.20)]{mrw}, and obtain $d_0(v_3^{(p^2-e(2))p^2+p+1})\cg v_3^{(p^2-e(2))p^2}(v_2z^p-v_2\pt2z) -v_2^{a_2}v_3^{(p^2-e(2))p^2-p}t_1^p(v_3^{p+1}+v_2z^p)+v_2^{p^2+p}v_3^{(p^2-e(2))p^2}t_2$ mod $J_{e(3)}$. The congruence on $d_0(x_4)$ follows from this and the congruence on $d_0(x_3)$ together with the definition of the element $x_3$. 

Inductively suppose that 
$$
d_0(x_{2k})\cg
v_2^{a_{2k}}x_{2k-1}^{p-1}t_1^p +v_2^{e(3)p^{2k-2}-e(2)}v_3^{(p^2-e(2))p^{2k-2}}(z-v_2^pt_2)\mod J_{e(3)p^{2k-2}}.
$$
Then, we calculate
\AR l{
d_0(x_{2k}^p)\cg \U{v_2^{pa_{2k}}x_{2k-1}^{(p-1)p}\tp2}a +v_2^{e(3)p^{2k-1}-e(2)p}v_3^{(p^2-e(2))p^{2k-1}}(\U{z^p}b-\U{v_2\pt2t_2^p}c)\\
d_0(-v_2^{pa_{2k}-1}x_{2k-1}^{(p-1)p}v_3) \\[-1mm]
\ak \uss{\kko{deriv}\\\kko{t1eRv3}}\cg -v_2^{pa_{2k}-1}x_{2k-1}^{(p-1)p}(\U{v_2\tp2}a-v_2^pt_1) +v_2^{e(3)p^{2k-1}-p-1}x_{2k-1}^{p^2-p-1}(z+\U{v_2t_2^p}c-v_2^pct_2) \\[-1mm]
d_0(-v_2^{e(3)p^{2k-1}-e(3)}v_3^{(p^2-e(2))p^{2k-1}+p+1})\cg -v_2^{e(3)p^{2k-1}-e(3)}v_3^{(p^2-e(2))p^{2k-1}}(\U{v_2z^p}b-v_2\pt2z)\\
\therefore\ d_0(x_{2k+1})\cg v_2^{pa_{2k}+p-1}x_{2k-1}^{(p-1)p}t_1+v_2^{e(3)p^{2k-1}-e(2)}v_3^{(p^2-e(2))p^{2k-1}}(2z-v_2^pct_2)
\qand\\
}
\AR l{
d_0(x_{2k+1}^p)\cg v_2^{pa_{2k+1}}x_{2k-1}^{(p-1)p^2}t_1^p+v_2^{e(3)p^{2k}-e(2)p}v_3^{(p^2-e(2))p^{2k}}(2z^p-v_2\pt2ct_2^p)\\
\ak \cg v_2^{pa_{2k+1}}(x_{2k+1}^{p-1}-v_2^{pa_{2k}-1}x_{2k-1}^{(p^2-p-1)p}v_3)t_1^p+v_2^{e(3)p^{2k}-e(2)p}v_3^{(p^2-e(2))p^{2k}}(2z^p-v_2\pt2ct_2^p)\\
\ak\cg v_2^{pa_{2k+1}}x_{2k+1}^{p-1}t_1^p+v_2^{e(3)p^{2k}-e(2)p}v_3^{(p^2-e(2))p^{2k}}(2z^p-v_2^{p^2-1}z-v_2^{p^2+p-1}t_2)\\
d_0(-2v_2^{e(3)p^{2k}-e(3)}v_3^{(p^2-e(2))p^{2k}+p+1})\cg -2v_2^{e(3)p^{2k}-e(3)}v_3^{(p^2-e(2))p^{2k}}(v_2z^p-v_2\pt2z) \\
}
\AR l{
\therefore\ d_0(x_{2k+2})\cg v_2^{pa_{2k+1}}x_{2k+1}^{p-1}t_1^p+v_2^{e(3)p^{2k}-e(2)}v_3^{(p^2-e(2))p^{2k}}(z-v_2^pt_2).
}
These complete the induction.

Put $d_0(x_i)\cg v_2^{a_i}(x_{i-1}^{p-1}\tp{\ep_i}+B_i+v_2^{a_{i-1}+1}C)$ mod $J_{e(3)p^{i-1}+1}$ for a cochain $C$.
It is easy to see $d_1(v_2^{a_i}(x_{i-1}^{p-1}\tp{\ep_i}+B_i))\cg 0$  mod $J_{e(3)p^{i-1}+1}$. It follows that $C$ is a cocycle of $\Om^1M^0_3$, and so $C$ represents a cohomology class $av_3^{c(i)}\ze_3\in H^1M^0_3$ for some $a\in \Z/p$ by \kko{HM}.
\q 

Put 
$$
d_0(x_i)\cg v_2^{a_i}A_i+v_2^{a_i}B_i\mbx{for $A_i=x_{i-1}^{p-1}\tp{\ep_i}$.}
$$
($\ep_i=\frac{1+(-1)^i}2$).
We introduce elements $y_i$ and $y_i'\in \Om^1E(3)_*$ by
$$
y_{s,i}= x_i^s\tp{\ep_{i+1}}-sx_i^{s-p+1}B_{i+1}, \qand
y_{s,i}'= x_i^s\tp{\ep_{i}}+\frac s2 v_2^{a_i}x_i^{s-1}A_i\tp{\ep_i} 
$$

\lem{yi}{For the elements $y_i$ and $y_i'$,
\sk{
d_1(y_{s,0})\cgs s(s+1)v_2^{2}v_3^{s-p-1}G_2\\
d_1(y_{s,1})\cgs s(s+1)v_2^{2p}v_3^{sp-2}G_0\\
d_1(y_{s,2})\cgs -s(s+1)v_2^{2p^2-p}v_3^{sp^2-2p}(t_1^p \ox z-v_2^{p}x) 
\\
d_1(y_{s,i})\cgs \cass{
-s(s+1)v_2^{2a_{2k+1}-p}x_{2k}^{sp-2}(t_1\ox z-v_2^{p}G_0)
&i=2k+1\\
-s(s+1)v_2^{2a_{2k+2}-p}x_{2k+1}^{sp-2}(2t_1^p\ox z-v_2^{p}K_0')
&i=2k+2,\\
}\qand\\
d_1(y_{s,1}')\cgs -sv_2^{p+1}v_3^{sp-2p}K_2\\
d_1(y_{s,2}')\cgs -sv_2^{p^2+p}v_3^{sp^2-p-1}K_0\\
d_1(y_{s,3}')\cgs sv_2^{a_3+p^2-p}v_3^{sp^3-p^2-p}(z\ox t_1 -v_2^p x')\\
d_1(y_{s,i}')\cgs \cass{
sv_2^{e(3)p^{i-2}-p-1}v_3^{(sp^2-p-1)p^{2k-2}}(z\ox t_1^p -v_2^pK_0)&i=2k\\
sv_2^{e(3)p^{i-2}-p-1}v_3^{(sp^2-p-1)p^{2k-1}}(2z\ox t_1 - v_2^pG'_0)&i=2k+1.
}
}
Here, $x=(t_2+t_1^{p+1})\ox t_1^p+t_1^p\ox t_1^{p+1}+\frac12t_1^{2p}\ox t_1$ and $x'=t_1^{p+1}\ox t_1+\frac12 t_1^p\ox t_1^{2}$, and
these congruences are considered modulo $J_{a+1}$, where $a$ is the largest power of $v_2$ in each congrucence.  
Furthermore, replace $K_0'$ and $K_0$ in the congruences on $d_1(y_{s,2k+2})$ and $d_1(y'_{s,2k})$
by $K_0'+v_2t_1^p\ox Z'$ and $K_0+v_2Z'\ox t_1^p)$, respectively, if we consider the congrucences modulo $J_{a+2}$.
}

\p
We note that $d_1(B_{i+1})\cg -d_1(A_{i+1})\cg -d_0(x_{i}^{p-1})\ox \tp{\ep_{i+1}}$ mod $I_2$ and $d_0(x_i^s)+sx_i^{s+1-p}d_0(x_{i}^{p-1})\cg \C {s+1}2 x_i^{s-2}d_0(x_i)^2$ mod $J_{3a_i}$. Indeed, $d_0(x_i^s)\cg sx_i^{s-1}d_0(x_i)+\C s2 x_i^{s-2}d_0(x_i)^2 $ mod $J_{3a_i}$.
We also see that $d_1(A_i\tp{\ep_i})\cg d_0(x_{i-1}^{p-1})\ox t_1^{2p^{\ep_i}}-2x_{i-1}^{p-1}\tp{\ep_i}\ox \tp{\ep_i}\cg d_0(x_{i-1}^{p-1})\ox t_1^{2p^{\ep_i}}-2A_i\ox \tp{\ep_i}$ mod $J_{a_{i-1}+2}$.
Then, we calculate 
\sk{
d_1(y_{s,i})\ukcgs{deriv} d_0(x_i^s)\ox \tp{\ep_{i+1}}-sd_0(x_i^{s+1-p})\ox B_{i+1}+sx_i^{s+1-p}d_0(x_{i}^{p-1})\ox \tp{\ep_{i+1}}\\
&\hsp{-.31in} \cg \C {s+1}2 x_i^{s-2}d_0(x_i)^2\ox \tp{\ep_{i+1}}-s(s+1)x_i^{s-p}d_0(x_i)\ox B_{i+1}\mod J_{2a_i+p}\\
d_1(y_{s,i}')\ukcgs{deriv} sx_i^{s-1}d_0(x_i)\ox \tp{\ep_{i}}+\frac s2v_2^{a_i}x_i^{s-1}d_0(x_{i-1}^{p-1})\ox t_1^{2p^{\ep_i}}-sv_2^{a_i}x_i^{s-1}A_i\ox \tp{\ep_i}\\
&\hsp{-.31in} \cg  sv_2^{a_i}x_i^{s-1}(B_i\ox \tp{\ep_{i}}+\frac12 d_0(x_{i-1}^{p-1})\ox t_1^{2p^{\ep_i}})  \mod J_{e(3)p^{i-2}+1}
}
Now we obtain the lemma from Lemma \ref{xi}. 
\q

\section{Proof of Lemma \ref{key}}
In this section, we define the cochains $\kt{p^i}s$ and verify the $d_1$-differential of them.

\subsection{The cochains $\kt {}{sp^{2k}}$ and $\kt p{sp^{2k+1}}$ for $s\in \Z_0$}
We define the cochains by
\sk{
\kt{}{s}\eqs y_{s,0}, \hsp{.5in} \kt{p}{sp}= y_{s,1},\\
\kt{}{sp^2}\eqs y_{s,2}-s(s+1)v_2^{2p^2-p}v_3^{sp^2-2p}\om_2,\\
\kt{p}{sp^{2k+1}}\eqs y_{s,2k+1}-s(s+1)v_2^{2a_{2k+1}-p}x_{2k}^{sp-2}\om_1\\
\kt{}{sp^{2k+2}}\eqs y_{s,2k+2}-s(s+1)v_2^{2a_{2k+2}-p^2-p}x_{2k+1}^{sp-2}(2\wt\om_2+v_2^{p^2}(2zt_1^p+v_2^{p}\dk_0))
}
for $k\ge 1$.
Then, the lemma for this case follows immediately from  Lemmas \ref{yi}, \ref{xi} and \ref{xz} toghere with \kko{homo}. Note also $2a_{2k+1}-p+2=2pa_{2k}+p$.

\subsection{The cochains $\kt {}{sp^{2k}}$ and $\kt p{sp^{2k+1}}$ for $s\in \Z_1$}
We put $s=tp^2-1$,
and define
 the cochains $\kt{}{(tp^2-1)p^{2k}}$ and $\kt{p}{(tp^2-1)p^{2k+1}}$ by
\sk{
v_2^{a_{2k+1}}\kt{}{(tp^2-1)p^{2k}}\eqs  -v_3^{(t-1)p^{2k+2}}w\pt{2k+1}-d_0(v_2^{p^{2k+1}-p^{2k-2}}v_3^{(tp^2-1)p^{2k}}\si_k)\\
&\hsp{.5in}+v_2^{p^{2k+2}-p^{2k-1}}v_3^{(tp-1)p^{2k+1}}W_k, \qand\\
%
\kt{p}{(tp^2-1)p^{2k+1}}\eqs \kt{}{(tp^2-1)p^{2k}}^p
}


\noindent
for the elements $\si_k$ in Lemma \ref{si}, $w$ in \kko{w13} and $W_k$ in Corollary \ref b.
Then, this case follows from
Lemmas \ref{si} and \ref w, Corollary \ref b and \kko{inta}.
We also use relations $w^{p^{2k+1}} 
\cg -v_2\pt{2k+1}v_3^{p^{2k+2}-p^{2k}}\tp{2k}$ and 
$
\b1{2k+3}\cg v_3^{(p-1)p^{2k+1}}\b1{2k}\mod I_3
$ given by \kko{t3=t}.

\subsection{The cochains $\kt {}{sp^{2k+1}}$ and $\kt p{sp^{2k}}$ for $s\in \Zpp$}

We begin with defining

\vsp{-.15in}

$$
\kt p{s}=v_3^st_1^p+sv_2v_3^{s-1}ct_2^p-s(s-1)v_2^2v_3^{s-2}\dk_1.
$$


\noindent
Then, we calculate by 
\kko{deriv}, \kko L, \kko{Det} and \kko{dgk}, and obtain 

\vsp{-.15in}

$$
d_1(\kt p{s})\cg s(s-1)v_2^2v_3^{s-2}K_1\mod J_3.
$$
%
%
Now we consider the cases for $p\mid s(s-1)$.

\vsp{-.05in}

\subsubsection{The cochains $\kt{p}{tp^k+1}$ for $k\ge 1$}
We define the cochains by
\sk{
\kt p{tp+1}\eqs v_3^{tp}z+tv_2^pv_3^{tp}t_2-tv_2^{p+1}v_3^{tp-p}ct_3^p,\\ 
\kt p{tp^{2}+1}\eqs x_{2}^{t}z+tv_2^{a_{2}}v_3^{(tp-1)p}\om_2,\\ 
\kt p{tp^{2k+1}+1}\eqs x_{2k+1}^{t}z+tv_2^{a_{2k+1}}v_3^{(tp-1)p^{2k}}\om_1
+tv_2^{a_{2k}+p+1}\kt p{(tp^2-1)p^{2k-1}}\qand\\ 
\kt p{tp^{2k+2}+1}\eqs x_{2k+2}^{t}z+tv_2^{a_{2k+2}-p^2}v_3^{(tp-1)p^{2k+1}}(\wt\om_2+v_2\pt2 zt_1^p)
}
%
in $\Om^1E(3)_*$ for $k\ge 1$, $t\in\Zpp$, $x_n$ in \kko{xi}, $z$ in \kko{At} 
 and $\om_i$ in Lemma \ref{xz}.
We verify this case by  a routine calculation using
 \kko{deriv}, \kko{L}, \kko{At}, \kko{Det} and \kko{Dec}.
We see that $\tp3\ox z\cg \eR(v_3)\tp3\ox t_1^p+v_2\tp3\ox ct_2^p-v_2^pv_3^{p-1}t_1\ox t_2$  and $\eR(v_3)\tp3\cg
v_3^pt_1+v_2ct_2\pt2$  mod $J_{p+1}$ by \kko{At}, \kko L and \kko{t3=t}.
It follows that $\tp3\ox z\cg -d_1(v_3^pt_2)+v_2d_1(ct_3^p)$ mod $J_{p+1}$, and then $d_1(v_3^{tp}z)\cg tv_2^pv_3^{tp-p}( -d_1(v_3^pt_2)+v_2d_1(ct_3^p))+\C t2v_2^{2p}v_3^{tp-1}t_1^{2}\ox t_1^p$ mod $J_{2p+1}$.
Thus, we obtain $d_1(\kt p{tp+1})$.

The congruences on $d_1(\kt p{tp^k+1})$ for $k\ge 2$ follow directly from Lemmas \ref{xi} and \ref{xz} and the results on $d_1(\kt p{(tp^2-1)p^{2k-1}})$ shown in the previous subsection. 

\vsp{-.05in}

\subsubsection{The cochains $\kt{p}{tp^k+e(k)}$ for $k\ge 2$}
We put $r=2n-1+\ep$ ($\ep\in\{0,1\}$), and 

\vsp{-.1in}

$$
\kt p{tp^r+e(r)}'=x_r^t\LR{w_{r+1}+v_2^{p^r-p^{r-3}}w_r\eR(\si_{n-1}\pt\ep)+v_2^{a_r}w_rt_1\pt\ep}
$$
for $w_r$ in \kko{wn}. Note that $w_r\cg v_3^{pe(r-2)}w_2\cg -v_3^{pe(r-2)}z$ mod $J_p$ by \kko{wn} and \kko{At}.
Then, $\kt p{tp^r+e(r)}'\cg x_r^tw_{r+1}\cg   -v_3^{tp^r+e(r)}t_1^p$ mod $I_3$. Furthermore, we calculate
to obtain

\vsp{-.15in}

$$
d_1(\kt p{tp^r+e(r)}')\cg -(t-1)v_2^{a_{r}}v_3^{tp^r+pe(r-2)}\tp\ep\ox z\mod J_{a_{r}+p}
$$
by Corollary \ref{3.7}, Lemmas \ref{si} and \ref{xi} together with \kko{deriv} and  \kko{inta}.
This case
now follows from Lemma \ref{xz} by setting $\kt p{tp^r+e(r)}=-\kt p{tp^r+e(r)}'+(t-1)$ $v_2^{a_{r}}v_3^{tp^r+pe(r-2)}\om_{1+\ep}$.
\q


\subsubsection{The cochains $\kt{p}{sp^{2k}}$ for $k\ge 1$  and $\kt{}{sp^{2k+1}}$ for $k\ge 0$}
We define $\kt{\ep_{i}}{sp^i}$ by

\vsp{-.2in}
\sk{
\kt{}{sp}\eqs y'_{s,1}, \hsp{.5in} \kt{p}{sp^2}= y'_{s,2},\\
%
\kt{}{sp^3}\eqs y'_{s,3}+sv_2^{e(3)p-p-1}v_3^{(sp^2-p-1)p}(zt_1-\om_1), \\
\kt{p}{sp^{4}}\eqs y'_{s,4}-\frac s2v_2^{e(3)p^{2}-p^2-p-1}v_3^{(sp^2-p-1)p^{2}}(\wt\om_2'-v_2\pt2z\tp{}), \\
\kt{}{sp^{2k+1}}\eqs y'_{s,2k+1}+2sv_2^{e(3)p^{2k-1}-p-1}v_3^{(sp^2-p-1)p^{2k-1}}(zt_1-\om_1),\qand \\
\kt{p}{sp^{2k+2}}\eqs y'_{s,2k+2}-sv_2^{e(3)p^{2k}-p^2-p-1}v_3^{(sp^2-p-1)p^{2k}}\wt\om_2.
}
where $\wt\om_2'=\wt\om_2-v_2^{p^2+p}t_1^pt_2-v_2^{e(3)}v_3^{-p^2}ct_4^p$.
Except for $d_1(\kt{p}{sp^{4}})$,
 the lemma for this case follows from Lemmas \ref{yi},  \ref{xz} with \kko{deriv}. 
 For $d_1(\kt{p}{sp^{4}})$, we use the relation
 $\tp4\ox \wt\om_2\cg v_3\pt2t_1^p\ox z-2v_2^pv_3\pt2K_0-v_2^{p+1}v_3\pt2\b10-d_1(v_2^pv_3\pt2t_1^pt_2+v_2^{p+1}ct_4^p)$ mod $J_{p+2}$  shown by an easy but tedious calculation using \kko{t3=t},
 \kko{2K0} and
 $\wt\om_2\cg -w_3\cg v_3^{p}(z+v_2^pt_2-v_2^{p+1}ct_3^p)$ mod $J_{p+2}$.
%

\subsection{The cochains $\kt {p^2}{tp-1}$ for $t\in \Z$}
Put 
%
%
$$
\kt {p^2}{tp-1}=-v_2^{-1}v_3^{(t-1)p}w.
$$
%
Then, the lemma for this case follows from Lemma \ref w.


\bibliographystyle{amsplain}

\end{document}